\documentclass[11pt]{article}

       \font\tenmsb=msbm10
       \font\sevenmsb=msbm7
       \font\fivemsb=msbm5
       \catcode`\@=11
       \ifx\amstexloaded@\relax\catcode`\@=\active
       \endinput\else\let\amstexloaded@\relax\fi
       \def\spaces@{\space\space\space\space\space}
       \def\spaces@@{\spaces@\spaces@\spaces@\spaces@\spaces@}
       \def\space@.{\futurelet\space@\relax}
       \space@. %
       \def\Err@#1{\errhelp\defaulthelp@\errmessage{AmS-TeX error: #1}}
       \def\relaxnext@{\let\next\relax}
       \def\accentfam@{7}
       
       \def\noaccents@{\def\accentfam@{0}}
       \def\Cal{\relaxnext@\ifmmode\let\next\Cal@\else
       \def\next{\Err@{Use \string\Cal\space only in math mode}}\fi\next}
       \def\Cal@#1{{\Cal@@{#1}}}
       \def\Cal@@#1{\noaccents@\fam\tw@#1}

       \def\Bbb{\relaxnext@\ifmmode\let\next\Bbb@\else
       \def\next{\Err@{Use \string\Bbb\space only in math mode}}\fi\next}
       
       \def\Bbb@#1{{\Bbb@@{#1}}}
       \def\Bbb@@#1{\noaccents@\fam\msbfam#1}
       \newfam\msbfam
       \textfont\msbfam=\tenmsb
       \scriptfont\msbfam=\sevenmsb
       \scriptscriptfont\msbfam=\fivemsb
\usepackage{hyperref}
\usepackage{tikz}
\usepackage{pgflibraryarrows}
\usepackage{pgflibrarysnakes}

\usepackage{fancyhdr}
\def\Z{{\Bbb Z}}
\def\R{{\Bbb R}}
\def\T{{\Bbb T}}

\newtheorem{thm}{Theorem}

\newtheorem{prop}{Proposition}

\newtheorem{rmk}{Remark}[section]

\newtheorem{notation*}{Notation}

\@addtoreset{equation}{section}

\newcommand{\qed}{\nolinebreak\hfill\rule{2mm}{2mm}
\par\medbreak}
\newcommand{\proof}{\par\medbreak\it Proof: \rm}

\newcommand{\beq}{\begin{equation} }
\newcommand{\eeq}{\end{equation} }

\begin{document}
\title{Mixed spectral types for the one frequency discrete quasi-periodic Schr\"{o}dinger operator}
\author{Shiwen Zhang\\
{\footnotesize Department of Mathematics, University of California Irvine}\\
{\footnotesize Irvine, CA 92617, U.S.A}\\
{\footnotesize Email: shiwez1@uci.edu}
\thanks{Research partially supported by NSF DMS-1401204.}}

\maketitle

\begin{abstract}
We consider a family of one frequency discrete analytic quasi-periodic Schr\"{o}dinger operators which appear in \cite{Bjer}. We show that this family provides an example of coexistence of absolutely continuous and point spectrum for some parameters as well as coexistence of absolutely continuous and singular continuous spectrum for some other parameters.
\end{abstract}


\section{Introduction}
In \cite{Bjer},  Bjerkl\"{o}v considers the following discrete quasi-periodic Schr\"{o}dinger operator on $l^2(\Z)$
       \begin{equation}\label{B-eg}
        (H_{K,\theta,\omega}u)_n=-u_{n+1}-u_{n-1}+V(\theta+n\omega)u_{n},\ \ n\in\Z
       \end{equation}
      where
       \begin{equation}\label{B-eg-V}
        V(\theta)=exp\Big(Kf(\theta+\omega)\Big)+exp\Big(-Kf(\theta)\Big)
       \end{equation}
      $\theta,\omega\in\T^b$, $f:\T^b\to\R$, is assumed to be a non-constant real-analytic
      function with zero mean, $\int_{\T^b}f(\theta){\rm d}\theta=0$ and $K\in\R$ is any constant. Consider the Lyapunov exponent $L(E)$ (see next section). In this explicit example, Bjerkl\"{o}v shows that for large $K$ we have a situation with mixed dynamics: zero Lyapunov exponent in a region close to $E = 0$ and positive for larger $E$.

      In this paper, we are going to show that for one frequency case, in Bjerkl\"{o}v's example (\ref{B-eg}), mixed dynamics actually lead to mixed spectra:  for some parameters $(\theta,\omega)$, $H_{K,\theta,\omega}$ has mixed absolutely continuous and point spectrum, and for some other $(\theta,\omega)$, $H_{K,\theta,\omega}$ has mixed absolutely continuous and singular continuous spectrum.

       Without loss of generality, we assume that $\|f\|_{C^1(\T)}=1$ and $f$ has analytic extension to the strip $|Im z|<h$, where $h>>{K}$ (e.g., $f$ can be taken as any entire function). It follows from \cite{Bjer} that $\min\{E\in\sigma(H_{K,\theta,\omega})\}=0$ for any $K,\theta,\omega$. And also it is not hard to show that $\max\{E\in\sigma(H_{K,\theta,\omega})\}\asymp e^{K\|f\|_\infty}$. For any $\epsilon>0$ small (w.l.g. we assume $0<\epsilon<1$), denote $I_{\epsilon,K}=[\epsilon, 4e^{K\|f\|_\infty}]$. We have $I_{\epsilon,K}\bigcap \sigma(H_{K,\theta,\omega})\neq\emptyset$. We say the frequency $\omega\in\T$ satisfies the Diophantine Condition  (denote by $\omega\in DC(\kappa,\tau)$) if
       $$\|\omega\cdot n\|\ge \frac{\kappa}{|n|^\tau},\ \  \forall n\in\Z\backslash\{0\}$$
for some $\kappa>0,\tau>0$. We say $\omega$ satisfies the Strong Diophantine Condition (denote by $\omega\in SDC(\kappa)$) if $$\|\omega\cdot n\|\ge \frac{\kappa}{|n|(\log(|n|+1))^2},\ \  \forall n\in\Z\backslash\{0\}$$
for some $\kappa>0$. Denote $DC=\bigcup_\kappa DC(\kappa,\tau)$ for some fixed $\tau>1$ and $SDC=\bigcup_\kappa SDC(\kappa)$. It is well known that $DC\supseteq SDC$ and both of them have full Lebesgue measure\footnote{
 Here $\|\cdot\|$ means the distance to the closest integer.
 }.

       The main results are as follows.
\begin{thm}
Let $V$ be given as in (\ref{B-eg-V}). Fix $\omega_0\in DC(\kappa,\tau)$. For any $\epsilon>0$, there are $K= K(\omega_0,\epsilon,f)>0$, $\delta=\delta(\omega_0,\epsilon,K)>0$, and for any $\omega\in B_{\delta}(\omega_0):=\{\omega\in\T:|\omega-\omega_0|<\delta\}$, there is $0<\epsilon_0=\epsilon_0(\omega,K,h,\|f\|_h)<\epsilon$ such that
 \begin{description}
   \item[(a)] for a.e. $\omega\in B_{\delta}(\omega_0)$ and a.e. $\theta\in\T$, $H_{K,\theta,\omega}$ has pure point spectrum in
         $I_{\epsilon,K}$ with exponentially decaying eigenvectors and has purely absolutely continuous spectrum in $[0,\epsilon_0]$.
   \item[(b)] for a.e. $\omega\in B_{\delta}(\omega_0)$, there is a dense $G_\delta$ set of $\theta$, such that $H_{K,\theta,\omega}$ has purely singular continuous spectrum in $I_{\epsilon,K}$ and has purely absolutely continuous spectrum in $[0,\epsilon_0]$.
 \item [(c)] for $\omega$ in a dense subset of $B_\delta(\omega_0)$ and for any $\theta$, $H_{K,\theta,\omega}$ has purely singular continuous spectrum in $I_{\epsilon,K}$ and has purely absolutely continuous spectrum in $[0,\epsilon_0]$.
         \end{description}
 \end{thm}

      Previously, Bourgain \cite{Bo} constructed quasi-periodic operator with two frequencies which has coexistence of absolutely continuous and point spectrum. While mixed spectra are expected to occur for generic one-frequency operators, such examples for the discrete case  have been considered difficult to construct explicitly.  Recently Bjerkl\"{o}v and Krikorian \cite{BK} announced an example of this nature. Avila in Theorem 13 \cite{A3} showed that in the neighborhood of the critical almost Mathieu operator, there are operators with at least $n$ alternances  between subcritical and supercritical regime for any $n$. Such operators are strong potential candidates for coexistence of p.p./s.c. and a.c. spectrum (with, moreover, many alternances). For continuous model, Fedotov and Klopp \cite{FK} showed coexistence of absolutely continuous and singular spectrum for a family of quasi-periodic operators and also gave a criterion for the existence of absolutely continuous and singular spectrum in the semi-classical regime.

      Here we give a short proof which shows that the operator (\ref{B-eg}) with potential (\ref{B-eg-V}) has mixed spectral types. The mixed nature of spectrum follows from coexistence of positive Lyapunov exponent and zero Lyapunov exponent which was obtained in \cite{Bjer} and a combination of several recent results on localization, reducibility and continuity .

\section{Singular spectrum in the positive Lyapunov exponent region}
Denote $$A(\theta,E)=
       \left(\begin{array}{cc}
              V(\theta)-E& -1 \\
              1 & 0
       \end{array}\right),\ \theta\in\T,\ \  A^n(\theta,E)=\prod^{0}_{k=n-1}A(\theta+k\omega,E),n>0$$
The Lyapunov exponent as usual (see \cite{CFKS}) is defined by
      $$ L(E)=\lim_{n\to \infty}\int_{\T}\frac{1}{n}\log\parallel A^n(\theta,E)\parallel{\rm d\theta}\ge 0$$

In the following, we would like to fix $f$ and consider the Lyapunov exponent $L(E,\omega,K)$ as function of energy $E$, frequency $\omega$, and parameter $K$. In \cite{Bjer}, Bjerkl\"{o}v proved that:
\begin{thm}[\cite{Bjer}]\label{Bjer-1}
        Assume that $V$ is as in (\ref{B-eg-V}), and that $\omega\in DC(\kappa,\tau)$.
        Then for any $\epsilon>0$ there is a $K_0 = K_0(\epsilon,f,\kappa,\tau)>0$ and $c=c(f)>0$
        such that for all $K >K_0$, we have
        $$L(E,\omega,K)\ge cK,\ \ \textrm{for all}\ \ E\not\in[0,\epsilon].$$
      \end{thm}

The proof is based on Large Deviation Theorem (LDT)-Avalanche Principle(AvP) scheme developed by Bourgain, Goldstein, Schlag \cite{BG,GS}. Due to some technical reasons, the largeness of  $K$ depends on the Diophantine Conditions of $\omega$ in this theorem, which means the positivity is not uniform for all $\omega$. However, we can get the following local non-perturbative positivity.  Bourgain and Jitomirskaya showed that Lyapunov exponent is jointly continuous in $(\omega,E)$ at any irrational frequency (Theorem 1, \cite{BJ}). The following result is obvious:

\begin{prop}\label{posi}
Fix any $\epsilon>0$ and $\omega_0\in DC(\kappa,\tau)$, let $K_0 = K_0(\epsilon,f,\omega_0)>0$
 be given as in Theorem \ref{Bjer-1}.  Then for any $K>K_0$, there is $\delta=\delta(\omega_0,\epsilon,K)>0$, such that for any $\omega\in B_{\delta}(\omega_0)$, $L(E,\omega,K)>0$ on $I_{\epsilon,K}$, where the lower bound only depends on $\omega_0,\epsilon,K,f$ and is uniform in $E$ and $\omega$.
\end{prop}

The absence of a.c. spectrum on $I_{\epsilon,K}$ is therefore obvious due to Kotani theory, see \cite{K} and discrete version in \cite{S}. What we want to claim is the pure point spectrum or purely singular continuous spectrum in this region.
\begin{description}
  \item[Anderson Localization (part (a))]  Let $\Omega=SDC\bigcap B_{\delta}(\omega_0)$, which is a full measure subset of $B_{\delta}(\omega_0)$. Notice that the positivity of $L(E,\omega,K)$ is uniform for $E\in I_{\epsilon,K}$ and $\omega\in \Omega$. Then according to the non-perturbative localization result of Bourgain and Goldstein (see Theorem 10.1 and Remark (3), Chapter 10 in \cite{Bo1}), for any $\theta\in \T$, a.e. $\omega\in \Omega$ , $H_{K,\theta,\omega}$ exhibits A.L. in
         $I_{\epsilon,K}$, i.e., $H_{K,\theta,\omega}$ has pure point spectrum restricted in $I_{\epsilon,K}$ and the corresponding eigenfunctions decay exponentially. Thus by Fubini's theorem, $H_{K,f,\theta,\omega}$ has A.L. for  a.e. $\omega\in \Omega$ and a.e. $\theta\in\T$.
  \item[Purely s.c. spectrum (part (b))] Let $\Omega$ be the same as in the previous part. Goldstein and Schlag \cite{GS2} show that for a.e. $\omega\in \Omega$, the intersection $\sigma(H_{K,f,\theta,\omega})\bigcap I_{\epsilon,K}$ is a Cantor set (see Theorem 1.1 in \cite{GS2}). Then according to a theorem of Gordon \cite{G2}, nowhere dense structure of the spectrum implies the absence of point spectrum for a dense $G_{\delta}$ set of $\theta$ (see Theorem 6 in \cite{G2}). Therefore, for a.e. $\omega\in \Omega$, there is a dense $G_{\delta}$ set of $\theta$ such that $H_{K,\theta,\omega}$ has purely singular continuous spectrum in $I_{\epsilon,K}$.
  \item[Purely s.c. spectrum (part (c))] Absence of point spectrum in this part is based on rational approximation.  More precisely, denote by $$\beta(\omega):=\limsup_n\frac{\log q_{n+1}}{q_n}$$ where $\frac{p_n}{q_n}$ is the $n^{th}$ rational approximation of $\omega$. Notice that
      $$\sup_n\sup_{(\theta,\omega)\in\T^2,\ E\in I_{\epsilon,K}}\frac{1}{|n|}\log\|A^n(\theta,E)\|\le 10K.$$
      Then by standard Gordon type argument (see e.g. \cite{G1,CFKS}), if $\beta(\omega)>40K$, then for any $\theta$, $H_{K,\theta,\omega}$ does not have any point spectrum. Combine with the positivity of Lyapunov exponent in $I_{\epsilon,K}$, the proof for purely s.c. spectrum of part (c) is completed. Notice that for any $\beta_0\in[0,\infty]$, the level set $\Omega_{\beta_0}:=\{\omega: \beta(\omega)=\beta_0\}$ is a dense set. For later purpose, we would like to pick the dense subset $\Omega_{\beta_0}\bigcap B_{\delta}(\omega_0)$ with $40K<\beta_0<h/2$.
\end{description}

\section{Absolutely continuous spectrum near the bottom}
        Next we are going to show that for any $\omega\in\T$ with finite $\beta(\omega)$, if the energy $E$ is sufficiently small (depends on $\omega$), then the Schr\"{o}dinger cocycle is almost reducible.  This will imply purely a.c. spectrum near the bottom of the spectrum for any phase $\theta$. To complete the proof of the main theorem, we first pick some $\omega$ near $\omega_0$ and some $\theta$ which give us point spectrum or singular continuous spectrum as in the previous part. Then for these pairs of $(\omega,\theta)$, we apply the almost reducibility result to get the coexistence of two types of spectrum.

        The key step to find purely a.c. spectrum near the bottom  is the following reducibility result at $E=0$, which generalizes Lemma 5.1 in \cite{Bjer} to the case $0<\beta(\omega)<\infty$.

         \begin{prop}\label{red}
         For any frequency $\omega$ with $\beta(\omega)<\infty$, if $h>2\beta$, then there exists analytic transformation $C:\T\to SL(2,\R)$ such that
        $$   C(\theta+\omega)A(\theta,0)C(\theta)^{-1}
           = A_0,
         $$
        where $$A_0= \left(
              \begin{array}{cc}
                1 & \hat{k} \\
                0 & 1 \\
              \end{array}
            \right),\ \ \hat{k}\in\R$$
         \end{prop}
         \begin{rmk}
         If $\omega$ is Diophantine or $\beta=0$, the proposition has been proved in Lemma 5.1 \cite{Bjer}. If $\beta>0$, we can still find such a transformation $C$ provided $h$ is large. The only loss is the decrease of the width of the analytic strip. Also the analytic norm of the transformation and the constant could be very large. Actually, $C$ has an analytic extension to the strip $|Im z|<h-2\beta$, with $\|C\|_{h-2\beta}\sim e^{K\|f\|_h}$. We also have $|\hat{k}|\sim e^{K\|f\|_h}$
         \end{rmk}
         \proof
         Recall the main steps in the proof of Lemma 5.1 \cite{Bjer},  if there are $g,h:\T\to\R$ satisfying the following equations
         \begin{eqnarray}
         % \nonumber to remove numbering (before each equation)
           g(\theta+\omega)-g(\theta) &=& f(\theta+\omega) \label{hom1}\\
           k(\theta)=-e^{-Kg(\theta-\omega)-Kg(\theta)} &,& \hat{k}=\int_{\T}k(\theta){\rm d}\theta \nonumber\\
           h(\theta+\omega)-h(\theta) &=& \hat{k}-k(\theta) \label{hom2}
         \end{eqnarray}
         then set
         $$C(\theta)=\left(
              \begin{array}{cc}
                1 & h(\theta) \\
                0 & 1 \\
              \end{array}
            \right)\cdot\left(
              \begin{array}{cc}
                0 & \exp\big(-Kg(\theta-\omega)\big) \\
                -\exp\big(Kg(\theta-\omega)\big) & \exp\big(Kg(\theta)\big) \\
              \end{array}
            \right)\in SL(2,\R)$$
        Direct computation shows that
         $$ C(\theta+\omega)A(\theta,0)C(\theta)^{-1}
           =  \left(
              \begin{array}{cc}
                1 & \hat{k} \\
                0 & 1 \\
              \end{array}
            \right) $$
         which is the desired form.\\

         For real analytic $f$ with zero average, if $\omega$ is Diophantine, equations (\ref{hom1}),(\ref{hom2}) always have real analytic solutions $g,h$, which is the case in \cite{Bjer}.

         If $\beta>0$, recall for $f$ analytic in the strip $|Im z|<h$, the Fourier coefficients of $f$ satisfy
         $|\hat{f}_k|\le \|f\|_he^{-h|k|}$,
         therefore, from Fourier series expansion, equation (\ref{hom1}) has an analytic solution $g$ in the strip $|Im z|<h-\beta$ provided $h>\beta$. From the definition of $k$, $k(\theta)$ also has analytic extension to the strip $|Im z|<h-\beta$ with
         $\|k\|_{h-\beta}\sim e^{K\|g\|_{h-\beta}}\sim e^{K\|f\|_{h}}$. Then for the same reason, equation (\ref{hom2}) also has an analytic solution $h$ in the strip $|Im z|<h-2\beta$ provided $h-\beta>\beta$. \qed

         Then it is easy to see that by applying $C$ to $A(\theta,E)$, we have
         \begin{eqnarray}
           % \nonumber to remove numbering (before each equation)
           C(\theta+\omega)A(\theta,E)C(\theta)^{-1}
           &=& A_0+C(\theta+\omega)\left(
              \begin{array}{cc}
                -E & 0 \\
                0 & 0 \\
              \end{array}
            \right)C(\theta)^{-1} \nonumber\\
            &=& A_0+EF(\theta) \in SL(2,\R) \label{almost1}
         \end{eqnarray}
         where
         $$F(\theta)=C(\theta+\omega)\left(
              \begin{array}{cc}
                -1 & 0 \\
                0 & 0 \\
              \end{array}
            \right)C(\theta)^{-1}.$$
            From the proof of Proposition \ref{red}, we see that $F$  has analytic extension to the strip $|Im z|<h-2\beta$ and  the largeness of $\|A_0\|,\|F\|_{h-2\beta}$ depend on $\omega,K,\|f\|_h$.\\

         For Diophantine frequency, as Bjerkl\"{o}v mentioned in Remark 2 in \cite{Bjer}, one can show purely absolutely continuous spectrum for sufficient small $E$ based on KAM approach as in \cite{E}. Here since we also need to deal with Liouvillean frequency, we want to prove all cases together with the  almost reducibility concept.\\

     We say the skew product system $(\omega,A)$ is \emph{almost reducible} if there exist $\eta>0$ and a sequence of analytic
maps $B^{(n)} : \T\to PSL(2;\R)$, admitting holomorphic extensions to the common
strip $|Im z|<\eta$ such that $B^{(n)}(z+\omega)A(z)B^{(n)}(z)^{-1}$ converges to a constant uniformly
in $|Im z|<\eta$. We need the following result about almost reducibility:
\begin{prop}[Corollary 1.2, \cite{A1}]\label{al}
Any one-frequency analytic quasi-periodic $SL(2, \R)$ cocycle close to
constant is analytically almost reducible.
\end{prop}
      \noindent   \emph{Proof of purely a.c. spectrum near the bottom}.
         According to Proposition \ref{al}, $(\omega,A_0+EF(\theta))$ is almost reducible for small $E$. More precisely, consider $A_0+EF$ in (\ref{almost1}). There exists $\epsilon_0=\epsilon_0(\omega,\|A_0\|,h,\|F\|_{h-2\beta})<\epsilon$ such that for $0<E<\epsilon_0$, $(\omega,A_0+EF(\theta))$ is almost reducible. (Such a quantitative version can be found in Theorem 1.2 \cite{HY} and Corollary 1.3 \cite{YZ}.) Therefore, $(\omega,A(\theta,E))$ is also almost reducible for $0<E<\epsilon_0$. As a corollary of almost reducibility \cite{A1,A2}
         \footnote{
         Corollary 1.6 in \cite{A1} established absolutely continuous spectrum as a consequence of almost reducibility for a.e. $\theta$, thus completing the proof of part a) and a.e. version of c) of Theorem 1. Part b) and extending part c) to all $\theta$ depend on the ``all $\theta$" version of this statement, to appear in \cite{A2}.    
         }, we have that for any $\theta$, $H_{K,f,\theta,\omega}$ has purely absolutely continuous spectrum in $[0,\epsilon_0]$. \qed

 ACKNOWLEDGEMENTS.
 The author
would like to thank Svetlana Jitomirskaya for many useful discussions. The author also
would like to thank the Isaac Newton Institute for Mathematical Sciences, Cambridge, for support and hospitality during the Programme \emph{Periodic and Ergodic spectral problems} where work on this paper was undertaken.

%%%%%%%%%%%%%%%%%%%%%%%%%%%%%%%%%%%%%%%%%%%%%%%%%%%%%%%%%%%%%%%%%%%%%%%%%%%%%%5
%%%%%%%%%%%%%%%%%%%%%%%%%%%%%%%%%%%%%%%%%%%%%%%%%%%%%%%%%%%%%%%%%%%%%%%%%%%%%%5
%%%%%%%%%%%%%%%%%%%%%%%%%%%%%%%%%%%%%%%%%%%%%%%%%%%%%%%%%%%%%%%%%%%%%%%%%%%%%%5

\end{document}